\documentclass[12pt]{amsart}
\usepackage{latexsym,amssymb,amscd,epic,eepic,epsfig,labelfig,psfig,graphics,color,ifthen,pstricks,verbatim,fancyhdr}
\usepackage{amsbsy,amsmath,amsfonts}
\usepackage{amsthm}
\usepackage[francais,english]{babel}
\usepackage{graphicx}
\usepackage[all]{xy}
\usepackage[latin1]{inputenc}

\newtheorem{Atheorem}{Theorem}

\newtheorem*{theoremnonum}{Theoreme}
\newtheorem{theorem}{Theorem}[section]
\newtheorem{lemma}[theorem]{Lemma}
\newtheorem{corollary}[theorem]{Corollary}
\newtheorem{proposition}[theorem]{Proposition}

\newtheorem*{probnonum}{Problem}
\newtheorem{definition}[theorem]{Definition}

\newtheorem{remark}[theorem]{Remark}

\newtheorem{NB}[theorem]{Nota Bene}

\newenvironment{sproof}[1]
{\begin{proof}[#1]} {\end{proof}}

\newcommand{\ov}[1]{\overline{#1}}
\newcommand{\N}{\mathbb N}

\newcommand{\Z}{\mathbb Z}

\newcommand{\Fm}{{\mathbb F}_m}

\newcommand{\Fi}{{\mathbb F}_{\infty}}
\newcommand{\Free}{{\mathbb F}}

\newcommand{\Mu}{{\mathcal M}}
\newcommand{\M}{{\mathcal M}_2}

\newcommand{\Cyc}{{\mathcal C}}

\newcommand{\Dihm}{{\mathcal D}_m}
\newcommand{\D}{{\mathcal D}}
\newcommand{\tD}{\tilde{\mathcal D}}

\newcommand{\Mm}{{\mathcal M}_m}

\newcommand{\V}[2]{V \left( #1 ,#2 \right)}
\newcommand{\Dn}{\mathbb{D}_{2n}}
\newcommand{\Dk}{\mathbb{D}_{2k}}
\newcommand{\Di}{\mathbb{D}_{\infty}}
\newcommand{\br}[1]{\lbrack #1 \rbrack}

\newcommand{\Pres}[2]{\left\langle{#1}\ \big\vert\ {#2}\right\rangle}

\newcommand{\Thu}{Th_{\forall}}

\newcommand{\U}{\mathfrak{U}}

\title{Limits of dihedral groups}
\author[Guyot]{Luc Guyot}
\address[L. Guyot]{Universit\'{e} de Gen\`{e}ve, Section de
Math\'{e}matiques, 2-4 rue du Li\`{e}vre, Case postale 64, 1211
Gen\`{e}ve 4, Switzerland} \email{Luc.Guyot@math.unige.ch}
\subjclass[2000]{Primary 20F05; Secondary 54G12, 11U09}
\begin{document}%%
%%%%%%%%%%%%%%%%%%
\thanks{This work was partially supported by the Swiss
National Science Foundation Grant $\sharp$~PP002-68627.}

\maketitle

\selectlanguage{english}

\begin{abstract}
We give a characterization of limits of dihedral groups in the space
of finitely generated marked groups. We also describe the
topological closure of dihedral groups in the space of marked groups
on a fixed number of generators.
\end{abstract}

\section{Introduction}
The space of marked groups is a topological setting for expressing
approximation among groups and algebraic limit processes in terms of
convergence. A convenient definition consists in topologizing the
set of normal subgroups of a given free group $\Free$ (we hence
topologize its marked quotients) with the topology induced by the
product topology of $\{0,1\}^{\Free}$. The idea of topologizing the
set of subgroups of a given group goes back to Chabauty's topology
on the closed subgroups of a locally compact group \cite{Chab50}. At
the end of its celebrated paper ''Polynomial growth and expanding
maps" \cite{Grom81}, Gromov sketched what could be a topology on
finitely generated groups and put it into practice in the purpose of
growth results. The space of marked groups on $m$ generators is
properly defined in \cite{Gri84} where the study of the neighborhood
of the first Grigorchuk group turned out to be fruitful. This
compact, totally disconnected and metrizable space has been used to
prove the existence of infinite groups with unexpected or rare
properties \cite{Ste96,Cham00,Sha00}. In \cite{ChGu05} Champetier
and Guirardel propose a new approach of Sela's limit groups
\cite{Sel01} (which appeared to coincide with the long-studied class
of finitely generated fully residually free groups, see definition
below) in the topological framework of the space of marked groups.
They give indeed the first proof not relying on the finite
presentability of limit groups \cite{KM98a,KM98b,Sel01,Gui04} that
limit groups are limit of free groups in the space of marked groups.
Among others, they provide a simple proof of the fact that a
finitely generated group is a limit group if and only if it has the
same universal theory as a free group of rank two. They hence relate
topology and logic. In the present paper, we use these links to
tackle the easier case of limits of dihedral groups. Our motivation
is:

\begin{probnonum}\rm\cite{Har00,ChGu05} \it
Describe the topological closure of finite groups in the space of
marked groups.
\end{probnonum}
We carry out such a description for the most elementary finite
groups: cyclic and dihedral finite groups.

Let $n$ be in $\{1,2,\dots\} \cup \{\infty \}$. We define the
dihedral group $$\Dn:=\Pres{a,b}{a^2=b^n=1,a^{-1}ba=b^{-1}},$$ where
we omit the relation $b^n=1$ when $n=\infty$. If $n$ is finite, then
$\Dn$ is a finite group of order $2n$. The groups $\mathbb{D}_2$
(which is isomorphic to $\Z/2\Z$) and $\mathbb{D}_4$ (which is
isomorphic to the Vierergruppe $\Z/2\Z \times \Z/2 \Z$) are the only
abelian dihedral groups. If $3 \le n<\infty$, then $\Dn$ is
isomorphic to the group of Euclidean isometries of a regular $n$-gon
$P_n$ (any function that maps $a$ to a reflection and $b$ to a
rotation with angle $\frac{2 \pi}{n}$, reflection and rotation both
preserving $P_n$, extends to a unique isomorphism from $\Dn$ to
$Isom(P_n)$). In this case we can identify $\{1,b,\dots,b^{n-1}\}$
with the group of rotations of $P_n$ and $\{a,ab,\dots,ab^{n-1}\}$
with the set of reflections of $P_n$. If $n \ge 3$ the center
$Z(\Dn)$ has two elements when $n$ is even ($1$ and $b^{n/2}$) and
$Z(\Dn)$ is trivial when $n$ is odd. For all $n<\infty$, the
subgroup $\langle b \rangle$ generated by $b$ is a cyclic normal
subgroup of order $n$ on which $\langle a \rangle$ acts by
conjugation. Thus $\Dn$ is isomorphic to the semidirect product
$\Z/n\Z \rtimes \Z / 2 \Z$ where the action of $\Z/2\Z$ on $\Z/n\Z$
is multiplying by $-1$. The infinite dihedral group $\Di$ is
centerless and is isomorphic to $\Z \rtimes \Z / 2 \Z$. All proper
quotients of $\Di$ are finite dihedral groups and any such group is
a quotient of $\Di$.

 We denote by $\Mu$ the space of all finitely
generated marked groups (see Section \ref{SecConvLogic} for a
definition). Unless otherwise stated, limits of groups considered
are limits in $\Mu$.

Let $P$ be a group theoretic property. A group $G$ is \it fully
residually \rm $P$ if for any finite subset $F$ of $G \setminus
\{1\}$ there is a group $H$ with $P$ and a homomorphism from $G$ to
$H$ that maps no element of $F$ to the trivial element. We address
the reader to Section \ref{SecConvLogic} for the definitions
required in our first theorem (convergence in $\Mu$, universal
theory $\Thu(G)$ of $G$, ultrafilter and ultraproduct). The only
abelian limits in $\Mu$ of dihedral groups are easily seen to be the
marked groups abstractly isomorphic to $\mathbb{D}_2$ or
$\mathbb{D}_4$ (see Section \ref{SecConvLogic}). We give in Section
\ref{SecChar} the following characterization of the non abelian
limits:
\begin{Atheorem}[Th. \ref{ThmLimDn}] Let $G$ be a non abelian finitely generated group. The following
conditions are equivalent:
\begin{itemize}
\item[$(\lim)$] $G$ is a limit of dihedral groups;
\item[$(res)$] $G$ is fully residually dihedral;
\item[$(iso)$] $G$ is isomorphic to a semidirect product
$A \rtimes \Z/2\Z$ where $A$ is a limit of cyclic groups on which
$\Z/2\Z$ acts by multiplication by $-1$;
\item[$(\Thu)$] $\Thu(G) \supset \bigcap_{n \ge 3} \Thu(\Dn)$;
\item[$(\Pi / \mathfrak{U})$] $G$
is isomorphic to a subgroup of $\left(\prod_{n \ge 3} \Dn \right)/
\mathfrak{U}$ for some ultra-filter $\mathfrak{U}$ on $\N$.
\end{itemize}
\end{Atheorem}
 Proposition \ref{PropLimZn} shows that limits of cyclic
groups are finitely generated abelian groups with cyclic torsion
subgroup.

We denote by $\omega$ the smallest infinite ordinal, \it{i.e.} \rm
the set $\N$ of positive integers endowed with its natural order. In
Section \ref{SecTop} we finally describe the set of limits of
dihedral groups on $m$ generators:
\begin{Atheorem}[Th.\ref{ThmTopDn}] \label{ThmB}
The topological closure $\Dihm$ of dihedral marked groups in $\Mm$
is homeomorphic to $\omega^{m-1}(2^m-1)+1$ endowed with the order
topology.
\end{Atheorem}
In other words, $\mathcal{D}_m$ is the disjoint union of $2^m-1$
copies of $\overline{\N}^{m-1}$ where $\overline{\N}=\N \cup
\{\infty\}$ is the Alexandroff compactification of $\N$. We use a
theorem of Mazurkiewicz and Sierpinski \cite{MS20} on
Cantor-Bendixson invariants of countable compact spaces to prove our
last result. For comparison, the set of abelian marked groups on $m$
generators is homeomorphic to $\omega^m+1$. This fact can be easily
derived from the proof of Theorem \ref{ThmB}.

\vspace{0.5cm}
\paragraph{\bf Acknowledgements. \rm}
I would like to thank Pierre de La Harpe, Goulnara Arjantseva and
Yves Stalder for their valuable comments on previous versions of
this article. I would also thank Thierry Coulbois and Vincent
Guirardel for interesting discussions and comments that originate
the subject of this article.

\section{Convergence and logic} \label{SecConvLogic}
Here, we give preliminary definitions and results that we need for
the proofs of the main theorems. We first define marked groups and
the topology (called Chabauty's topology, Cayley's topology or weak
topology) on the set of marked groups. Second, we relate convergence
in the space of marked groups to universal theory and ultraproducts
(see \cite{ChGu05} for a more full-bodied exposition).
\begin{definition}
The pair $(G,S)$ is a marked group on $m$ generators if $G$ is group
and $S=(g_1,\dots,g_m)$ is an ordered system of generators of $G$.
We call also $S$ a generating $G$-vector of length $n$ (or simply a
marking of $G$ of length $m$) and we denote by $\V{G}{m}$ the set of
these $G$-vectors.
\end{definition}
We denote by $r(G)$ the \it rank \rm of $G$, \it i.e. \rm the
smallest number of generators of $G$. Two marked groups $(G,S)$ and
$(G',S')$ (with $S=(g_1,\dots,g_m)$ and $S'=(g_1',\dots,g_m')$) are
\it equivalent \rm if there is an isomorphism $\phi:G
\longrightarrow G'$ such that $\phi(g_i)=g_i'$ for $i=1,\dots,m$.
\begin{NB}
We denote also by $(G,S)$ the equivalence class of $(G,S)$ and we
call this class a marked group.
\end{NB}

Let $\Fm$ be the free group with basis $(e_1,\dots,e_m)$ and let
$\Fi$ the free group with basis $(e_i)_{i \ge 1}$. Let $G$ and $G'$
be two groups and let $p:\Fm\twoheadrightarrow G$ and
$p':\Fm\twoheadrightarrow G'$ be two epimorphisms. The epimorphisms
$p$ and $p'$ are \it equivalent \rm if there is an isomorphism
$\phi$ such that the following diagram commutes: $\xymatrix{\Fm
\ar@{->>}[r]^{p} \ar@{->>}[rd]^{p'}& G \ar[d]^{\phi}
\\&G'}$

Using the universal property of $\Fm$, we establish a natural
bijection from marked groups on $m$ generators to equivalence
classes of epimorphisms with source $\Fm$. The last set is clearly
in one-one correspondence with the set $\mathcal{N}(\Fm)$ of all
normal subgroups of $\Fm$. We endow $\mathcal{N}(\Fm)$ with the
topology induced by the Tychonoff topology on the product
$\{0,1\}^{\Fm}$ which identifies with the set of all subsets of
$\Fm$. We easily check that $\mathcal{N}(\Fm)$ is a closed subspace
of $\{0,1\}^{\Fm}$. We hence define a compact, metrizable and
totally discontinuous topology on the corresponding set $\Mm$ of
marked groups on $m$ generators. We have for instance $(\Z/k\Z,
\br{1}_k) \underset{k \to \infty}{\longrightarrow} (\Z, 1)$ in
$\Mu_1$ and $(\Dk, (a,b)) \underset{k \to \infty}{\longrightarrow}
(\Di,(a,b))$ in $\M$ with respect to the Chabauty's topology. Let
$G$ be a finitely generated group. We define $\mathcal{N}(G)$ as the
set of normal subgroups of $G$ endowed with the topology induced by
$\{0,1\}^G$. Let $p:\Fm \twoheadrightarrow G$ be an epimorphism and
let $p^{\ast}:\mathcal{N}(G) \longrightarrow \mathcal{N}(\Fm)$ be
defined by $p^{\ast}(N)=p^{-1}(N)$. The map $p^{\ast}$ is clearly
continuous and injective. When $G$ is finitely presented, we have:
\begin{lemma}\rm\cite[Lemma 2.2]{ChGu05}\it \label{LemPresFinQuo}
The map $p^{\ast}$ is an open embedding.
\end{lemma}
Thus, for any marking $S$ of a finitely presented group $G$ in
$\Mm$, there is a neighborhood of $(G,S)$ containing only quotients
of $G$ (the quotient map is induced by the natural bijection between
the markings). Namely, this is the set corresponding to the image of
$p^{\ast}$ where $p$ is the epimorphism defined by means of $S$.
More generally, $p:G\twoheadrightarrow H$ induces an open embedding
$p^{\ast}:\mathcal{N}(H) \hookrightarrow \mathcal{N}(G)$ if and only
if $\ker p$ is finitely generated as a normal subgroup of $G$.
\begin{remark} \label{RemFinIso}
It follows from Lemma \ref{LemPresFinQuo} that finite groups are
isolated in $\Mu$. Isolated groups are characterized in \rm
\cite{Gri05,dCGP07}\it. Since nilpotent groups are finitely
presented, we can also derive of Lemma \ref{LemPresFinQuo} that the
set of nilpotent groups of a given class $c$ is open in $\Mu$.
\end{remark}
From the spaces $\Mm$, we build up the space $\Mu$ of all finitely
generated marked groups. Observe first that the map
$(G,(g_1,\dots,g_m)) \mapsto (G,(g_1,\dots,g_m,1))$ defines a
continuous and open embedding $i_n:\Mm \hookrightarrow \Mu_{m+1}$.
The inductive limit $\Mu$ of the of the system $\{i_n:\Mm
\hookrightarrow \Mu_{m+1}\}_{m \ge 1}$
 is a metrizable locally compact and totally discontinuous space.
 Observe that a convergence in $\Mu$ boils down to a convergence in
 $\Mm$ for some $m$ and the converse is obvious. The space $\Mu$
 can be viewed as the set of (equivalence classes) of groups marked with
 an infinite sequence of generators which are eventually trivial. Let
 $v_1,\dots,v_k,w_1,\dots,w_l$ be
elements of $\Fi$ and let $(\Sigma)$ be
the system $\left\{ \begin{array}{c} v_1 = 1,\dots,v_k = 1 \\
w_1 \neq 1,\dots ,w_l \neq 1 \end{array} \right.$. We denote by
$O_{\Sigma}$ the set of marked groups $(G,S)$ of $\Mu$ for which $S$
(possibly completed by trivial elements) is a solution of $(\Sigma)$
in $G$. The family of sets $O_{\Sigma}$ defines a countable basis of
open and closed subsets of $\Mu$.
\begin{remark} \label{RemRes}
Let $P$ be group theoretic property stable under taking subgroups
and let $G$ be a fully residually $P$ group.
 Using the previous basis, we easily deduce
that for any marking $S$ of $G$, $(G,S)$ is the limit in $\Mu$ of a
sequence of marked groups with $P$. The converse is also true when
$G$ is finitely presented because of Lemma \ref{LemPresFinQuo}.
\end{remark}

\begin{lemma}[Subgroup and convergence \cite{ChGu05}]
\label{LemSubConv} Let $(G_n,S_n)_n$ be convergent sequence in $\Mu$
with limit $(G,S)$. Let $H$ be a finitely generated subgroup of $G$.
Then, for any marking $T$ of $H$ there is a sequence of marked
groups $(H_n,T_n)_n$ which converges in $\Mu$ to $(H,T)$ and such
that $H_n$ is a subgroup of $G_n$.
\end{lemma}

It is important to note that being a limit in $\Mu$ of marked groups
with a given property $P$ does not depend on the marking:

\begin{lemma}\rm\cite{Cham00,dCGP07}\it \label{LemIsoHom}
Let $(G,S)$ and $(H,T)$ be in $\Mu$. Assume that $G$ and $H$ are
abstractly isomorphic. Then we can find a neighborhood $U$ of
$(G,S)$, a neighborhood $V$ of $(H,T)$, and a homeomorphism $\phi:U
\longrightarrow V$  mapping $(G,S)$ onto $(H,T)$. Moreover, $\phi$
preserves the isomorphism relation.
\end{lemma}
Thus, being a limit of finite (or equally free, cyclic, dihedral)
groups is a group property that doesn't depend on the marking. As
$\Di$ is a limit of finite dihedral groups, limits of dihedral
groups are all limits of finite ones.

Let $G$ be a finite group. For all $m \ge r(G)$, there are only
finitely many marked groups isomorphic to $G$ in $\Mm$ (there are
only finitely many equivalence classes of epimorphisms $p:\Fm
\twoheadrightarrow G$). Hence any compact subset of $\Mu$ contains
only finitely many marked groups isomorphic to $G$. Since the
property of being abelian is open in $\Mu$ (use Lemma
\ref{LemPresFinQuo} for instance), abelian limits of dihedral groups
are limits of (finite) abelian dihedral groups $\mathbb{D}_2$ or
$\mathbb{D}_4$. These limits are consequently isomorphic either to
$\mathbb{D}_2$ or $\mathbb{D}_4$. For this reason, we consider only
non abelian limits.

We turn now to relations between convergence and logic. We fix a
countable set of variables $\{x_1,x_2,\dots\}$ and the set of
symbols $\{\wedge,\vee,\neg,\forall,(,),\cdot,^{-1},=,1\}$ that
stand for the usual logic and group operations. We consider the set
$\Thu$ of \it universal sentences \rm in group theory written with
the variables, \it i.e. \rm all formulas $\forall x_1 \dots \forall
x_k \, \phi(x_1,\dots,x_k)$, where $\phi(x_1,\dots,x_k)$ is a
quantifier free formula built up from the variables and the
available symbols (see \cite{BeSl69,ChGu05} for a precise
definition). For instance, $\forall x \forall y \,( xy=yx)$ is a
universal sentence that is true in any abelian group while $\forall
x \, (x=1 \vee x^2 \neq 1)$ expresses that there is no $2$-torsion
in a group. If a universal sentence $\sigma$ is true in $G$ (we also
say that $G$ \it satisfies \rm $\sigma$), we write $G \models
\sigma$. We denote by $\Thu(G)$ the set of universal sentences which
are true in $G$. Let $(A_n)_n$ be a sequence of subsets of $\Thu$.
We set $\limsup A_n:=\bigcap_n \bigcup_{k \ge n}A_k$ and $\liminf
A_n:=\bigcup_n \bigcap_{k \ge n} A_k$. The set of marked groups of
$\Mu$ satisfying a given family of universal sentences define a
closed subset. More precisely, quoting Proposition 5.2 of
\cite{ChGu05} 5.3 and reformulating slightly Proposition 5.3, we
have:
\begin{proposition}[Limits and universal theory] \label{PropConvThu}
Let $G$ be a finitely generated group and let $(G_n)_n$ be a
sequence of finitely generated groups.
\begin{itemize}
\item[$(i)$] Assume that $(G,S)$ is the limit in $\Mu$ of $(G_n,S_n)_n$
for some ordered generating set $S$ of $G$ and some ordered
generating set $S_n$ of $G_n$. Then $\Thu(G) \supset \limsup
\Thu(G_n)$.
\item[$(ii)$]
Assume $\Thu(G) \supset \bigcap_n \Thu(G_n)$. Then, for any ordered
generating set $S$ of $G$, there is some integer sequence $(n_k)_k$
such that $(G,S)$ is the limit in $\Mu$ of some sequence $(H_k,T_k)$
satisfying $H_k \le G_{n_k}$.
\end{itemize}
\end{proposition}
It directly follows that a variety of groups (see \cite{Neu67} for a
definition) defines closed subspaces of $\Mm$ and $\Mu$. For
example, limits of dihedral groups are metabelian groups
($2$-solvable groups) because dihedral groups are metabelian. It can
be easily proved, by using Lemma \ref{LemPresFinQuo}, that a variety
of groups defines an open subspace of $\Mm$ if and only if its free
group on $m$ generators is finitely presented.

Convergence of a sequence $(G_n)_n$ in $\Mu$ can also be related to
ultraproducts of the $G_n$'s that we define now.

\begin{definition}[Ultrafilter]
An ultrafilter $\mathfrak{U}$ on $\N$ is a finitely additive measure
with total mass $1$ which takes values in $\{0,1\}$. In other words,
it is map from $\mathcal{P}(\N)$ (the set of all subsets of $\N$) to
$\{0,1\}$ satisfying $\U(\N)=1$ and such that, for all disjoint
subsets $A$ and $B$ of $\N$, we have $\U(A \cup B)=\U(A)+\U(B)$.
\end{definition}
 Let $\U$ be an ultrafilter on $\N$ and let $(G_n)_n$ be
sequence of groups. There is a natural relation on the cartesian
product $\prod_n G_n$: $(g_n)_n$ and $(h_n)_n$ are \it equal
$\U$-almost everywhere \rm if $\U(\{n \in \N \, \vert \, g_n=g'_n
\})=1$.
\begin{definition}[Ultraproduct]
The ultraproduct of the sequence $(G_n)_n$ relatively to $\U$ is the
quotient of $\prod_n G_n$ by the equivalence relation of equality
$\U$-almost everywhere. We denote it by $\left(\prod_n G_n
\right)/\U$.
\end{definition}

Let $(A_n)_n$ be a sequence of subsets of $\Thu$. We set $\lim_{\U}
A_n:=\{\phi \in \Thu \, \vert \, \phi \mbox{ belongs to } A_n \mbox{
 for }\U-\mbox{almost every }n\}$. We have:
\begin{theorem}[L\v{o}s \cite{BeSl69}]
$\Thu\left(\left(\prod_n G_n \right)/\U\right)=\lim_{\U}\Thu(G_n)$.
\end{theorem}
An ultrafilter is said to be \it principal \rm if it is a Dirac
mass. The following proposition relates convergence in $\Mu$ to
ultraproducts:
\begin{proposition}\rm (Limits and ultraproducts \cite[Prop. 6.4]{ChGu05}) \it  \label{PropUltra}
Let $(G,S)$ be in $\Mu$ and let $(G_n,S_n)_n$ be sequence in $\Mu$.
\begin{itemize}
\item[$(i)$] If $(G_n,S_n)_n$ accumulates on $(G,S)$
in $\Mu$, then $G$ embeds isomorphically into $\left(\prod_n
G_n\right)/\mathfrak{U}$ for some non principal ultrafilter $\U$.
\item[$(ii)$] If $(G_n,S_n)_n$ converges to $(G,S)$ in $\Mu$, then $G$
embeds isomorphically into $\left(\prod_n G_n\right)/\U$ for all
 non principal ultrafilter $\U$.
\item[$(iii)$] Let $H$ be a finitely generated group. If $H$
embeds isomorphically into $\left(\prod_n G_n\right)/\U$ for some
non principal ultrafilter $\U$, then for all ordered generating set
$T$ of $H$, we can find a sequence of integers $(n_k)_k$ and a
sequence $(H_k,T_k)_k$ that converges to $(H,T)$ in $\Mu$ and such
that $H_k \le G_{n_k}$ for all $k$.
\end{itemize}
\end{proposition}

\section{Cantor-Bendixson invariants} \label{SecCB}
This section is devoted to the basics of the Cantor-Bendixson
analysis we use in Theorem \ref{ThmTopDn}. Let $X$ be a topological
space. We denote by $X'$ the set of accumulation points of $X$. Let
$X^{(0)}:=X$. We define by transfinite induction the $\alpha$-th
derived set of $X$: $X^{(\alpha)}=(X^{\alpha-1})'$ if $\alpha$ is a
successor and $X^{(\alpha)}=\bigcap_{\beta<\alpha}X^{\beta}$ if
$\alpha$ is a limit ordinal. We denote by $\omega$ the set $\N$ of
integers endowed with its natural order. We use the following
topological classification theorem:
\begin{theoremnonum}[Mazurkiewicz-Sierpinski Theorem \cite{MS20}]
 For any given pair
$(\alpha,n)$ where $\alpha$ is countable ordinal number and $n$
belongs to $\N$, there is (up to homeomorphism) a unique countable
compact space $X$ such that $X^{(\alpha)}$ has exactly $n$ points:
the set $\omega^{\alpha}n+1$ endowed with the order topology.
\end{theoremnonum}
 The pair $(\alpha,n)$ is \it the characteristic system of \rm $X$. For example, the
Alexandroff compactification $\overline{\N}$ of $\N$ is homeomorphic
to $\omega+1$. Its characteristic system is then $(1,1)$. Similarly,
the caracteristic system of $\ov{\N}^k$ is $(k,1)$. The \it
Cantor-Bendixson rank of a point \rm $x$ in $X$ is the smallest
ordinal number $\alpha$ such that $x$ doesn't belong to
$X^{(\alpha)}$. A countable compact space has the characteristic
system $(\alpha,n)$ if and only if the set of points of maximal
Cantor-Bendixson rank (\it i.e. \rm rank $\alpha$) has cardinal $n$.
Observe that the Cantor-Bendixson rank of $(G,S)$ in $\Mm$ is the
Cantor-Bendixson rank of $(G,S)$ in $\Mu$ and does not depend on $S$
because of Lemma \ref{LemIsoHom}.

\section{Characterization of limits} \label{SecChar}
\begin{theorem}\label{ThmLimDn} Let $G$ be a non abelian finitely generated group. The following
conditions are equivalent:
\begin{itemize}
\item[$(\lim)$] $G$ is a limit of dihedral groups;
\item[$(res)$] $G$ is fully residually dihedral;
\item[$(iso)$] $G$ is isomorphic to a semi-direct product
$A \rtimes \Z/2\Z$ where $A$ is a limit of cyclic groups on which
$\Z/2\Z$ acts by multiplication by $-1$;
\item[$(\Thu)$] $\Thu(G) \supset \bigcap_{n \ge 3} \Thu(\Dn)$;
\item[$(\Pi / \mathfrak{U})$] $G$
is isomorphic to a  subgroup of $\left(\prod_{n \ge 3} \Dn \right)/
\mathfrak{U}$ for some ultra-filter $\mathfrak{U}$ on $\N$.
\end{itemize}
\end{theorem}
We first give a characterization of limits of cyclic groups that we
use in the proof of theorem \ref{ThmLimDn}:

\begin{proposition}\label{PropLimZn} Let $G$ be a finitely generated group. The following
conditions are equivalent:
\begin{itemize}
\item[$(\lim)_c$] $G$ is a limit of cyclic groups;
\item[$(res)_c$] $G$ is fully residually cyclic;
\item[$(iso)_c$] $G$ is isomorphic to an abelian group with cyclic
torsion subgroup;
\item[$(\Thu)_c$] $\Thu(G) \supset \bigcap_{n \ge 1} \Thu(\Z/n\Z)$;
\item[$(\Pi / \mathfrak{U})_c$]$G$ is isomorphic to a subgroup of $\left(\prod_{n \ge 1}
\Z/n\Z \right)/ \mathfrak{U}$ for some ultra-filter $\mathfrak{U}$
on $\N$.
\end{itemize}
\end{proposition}

We now show Proposition \ref{PropLimZn} and then Theorem
\ref{ThmLimDn}.

\begin{sproof}{Proof of Proposition \ref{PropLimZn}}
Here is the logical scheme of the proof: $$\xymatrix{
(res)_c\ar@{=>}[dr]&&(\Pi/\U)_c\ar@{=>}[dd]
\\&(\lim)_c
\ar@{=>}[ur]\ar@{=>}[dl]&\\(iso)_c\ar@{=>}[uu]&&(Th_{\forall})_c\ar@{=>}[ul]}$$

We begin with the right triangle of implications.

 $(\lim)_c
\Longrightarrow (\Pi/\U)_c$ : we first assume that there is a
sequence $(\Z/n\Z,S_n)_n$ in $\Mu$ which accumulates on $(G,S)$ for
some ordered generating set $S$ of $G$. Then $G$ embeds
isomorphically into $\left(\prod_{n \ge 1} \Z/n\Z
\right)/\mathfrak{U}$ for some non principal ultrafilter $\U$ by
Proposition \ref{PropUltra}$(i)$. If $G$ is the limit of stationary
sequence of finite cyclic groups, then $G$ is a finite cyclic group
isomorphic to $\Z/k\Z$ for some $k \ge 1$. We then set $\U$ as the
dirac mass in $k$.

$(\Pi/\U)_c \Longrightarrow (Th_{\forall})_c$: as $G$ embeds in
$\left(\prod_{n \ge 1} \Z/n\Z \right)/\mathfrak{U}$ for some
ultrafilter $\U$, we get then $Th_{\forall}(G) \supset
Th_{\forall}(\left(\prod_{n \ge 1} \Z/n\Z \right)/\mathfrak{U})$. By
L\u{o}s's theorem, we have $Th_{\forall}(\left(\prod_{n \ge 1}
\Z/n\Z \right)/\mathfrak{U})=\lim_{\U} Th_{\forall}(\Z/ n\Z)$.
\\Since $\lim_{\U} Th_{\forall}(\Z/ n\Z) \supset \bigcap_{n \ge 1}
Th_{\forall}(\Z/ n\Z)$, the result follows.

$(Th_{\forall})_c \Longrightarrow (\lim)_c$ : as $\Thu(G) \supset
\bigcap_{n \ge 1} \Thu(\Z /n\Z)$, we get by Proposition
\ref{PropConvThu}$(ii)$ that $G$ is a limit in $\Mu$ of subgroups of
cyclic groups. Hence $G$ is a limit of cyclic groups.

We now consider the left triangle.

$(\lim)_c\Longrightarrow (iso)_c$: as $G$ is a limit of cyclic
groups, $G$ is abelian (recall that the property of being abelian is
closed in $\Mu$). By the subgroup convergence lemma, the torsion
subgroup $Tor(G)$ of $G$ is itself a limit of cyclic groups. Because
$Tor(G)$ is finite, it is isolated in $\Mu$ (Rem. \ref{RemFinIso}).
Any sequence converging to $Tor(G)$ is then a stationary sequence.
As a result, $Tor(G)$ is cyclic.

$(iso)_c  \Longrightarrow (res)_c$ : we write $G=\Z^r \times \Z/k\Z$
with $r \ge 0$ and $k \ge 1$. Let $F$ be a finite subset of
$G\setminus \{0\}$ and let $E \subset \Z$ be the set of all $i$-th
coordinates of elements of $F$ for all $1 \le i \le r$. Choose
distinct primes $p_1,\dots,p_r$ such that each $p_j$ is coprime with
all elements of $E \cup \{k\}$. The quotient map $q :G
\twoheadrightarrow \left(\prod_{j=1}^r \Z/p_j\Z \right) \times \Z /
k \Z$ defined in an obvious way has cyclic image by the Chinese
Theorem. Clearly, no element of $F$ is mapped to the trivial element
by $q$.

$ (res)_c  \Longrightarrow (\lim)_c$: follows from Remark
\ref{RemRes}.
\end{sproof}

\begin{sproof}{Beginning of the proof of Theorem \ref{ThmLimDn}}

The first part of the proof is similar to the proof of Proposition
\ref{PropLimZn}. Actually, the following implications can be shown
in the very same way:
 $$\xymatrix{
(res)\ar@{=>}[dr]&&(\Pi/\U)\ar@{=>}[dd]
\\&(\lim)
\ar@{=>}[ur]&\\(iso)\ar@{=>}[uu]&&(Th_{\forall})\ar@{=>}[ul]}$$ The
only noticeable changes occure in the proof of
$(\Thu)\Longrightarrow(\lim)$: if $G$ is isomorphic to an abelian
dihedral group then $G$ is obviously the limit of a stationary
sequence of dihedral groups. We can assume then that $G$ is not
abelian. Since the property of being abelian is open in $\Mu$, $G$
is the limit of non abelian subgroups of dihedral groups. Hence, $G$
is the limit of dihedral groups.

 We then complete the
proof  by showing: $(Th_{\forall})\Longrightarrow (iso)$. This last
step relies on specific sentences that can be found in the universal
theory of all non abelian dihedral groups. We use the following
lemma:

\begin{lemma} \label{LemP}The following sentences are true in any non abelian
dihedral group:
\begin{itemize}

\item[$(P_1)$] $\forall x \forall y \, ( x^2 \neq 1 \wedge y^2 \neq 1)\Rightarrow
xy=yx$ (\it rotations commute\rm);

%\item[$(P_2)$] $\forall x \forall y \forall z \forall t \, (xz \neq zx \wedge yt \neq ty \wedge xy=yx)\Rightarrow
%(x^2=1 \wedge y^2=1) \vee (x^2 \neq 1 \wedge y^2 \neq 1)$ (\it two
%non central commuting symmetries are either both reflection or
%rotations\rm);

\item[$(P_2)$] $\forall x \forall y \forall z \, (x \neq 1 \wedge x^2 = 1 \wedge y^2 \neq 1 \wedge xz \neq zx )\Rightarrow
x^{-1}yx=y^{-1}$ (\it conjugation of a rotation by a reflection
reverses its angle\rm);

%\item[$(P_3)$] $\forall x \forall y \forall z \,( x^2 \neq 1 \wedge y^2 \neq 1   \wedge (xy)^2 = 1) \Rightarrow
%(xy)z=z(xy)$ (\it a rotation of order $2$ is central\rm);

\item[$(P_3)$] $\forall x \forall y \forall z \forall t \forall u \, (xz \neq zx \wedge yt \neq ty \wedge
 x^2=1 \wedge y^2=1 \wedge (xy)^2=1)\Rightarrow (xy)u=u(xy)$ (\it the product of two
 commuting reflections is central\rm);

\item[$(P_4)$] $\forall x \forall y \forall z \forall t \,(x \neq 1 \wedge x^2 = 1 \wedge y \neq 1 \wedge y^2 =1
 \wedge z^2 \neq 1 \wedge t^2 \neq 1 \wedge xz=zx \wedge yt=ty) \Rightarrow
x=y$ (\it there is at most one central element of order $2$ \rm).
\end{itemize}
\end{lemma}

\begin{sproof}{Proof of Lemma \ref{LemP}}
There are two kinds of symmetries in a dihedral group $\Dn$ ($n \ge
3$): the rotations (positive isometries of the Euclidean plane or
the Euclidean line) and the reflections (negative isometries). The
non central rotations $x$ of $\Dn$ are characterized by the
inequation $x^2 \neq 1$. All reflections have order $2$ and there is
possibly one central rotation of order $2$. All sentences can be
readily shown by writing $\Dn$ as the semidirect product $\Z/n\Z
\rtimes \Z/2\Z$
\end{sproof}

In fact, the three first sentences are true in any \it generalized
dihedral group \rm $Dih(A):=A \rtimes \Z/2\Z$ where $A$ is abelian
and $\Z/2\Z$ acts on $A$ by multiplication by $-1$.\\

\it End of the proof of Theorem \ref{ThmLimDn} \rm

By assumption, all sentences of Lemma \ref{LemP} are true in $G$.
 We denote by $A$ the subgroup of $G$ generated by the set $\{x \in G \,
\vert \, x^2 \neq 1\} \cup Z(G)$. We show the following claims:
\begin{itemize}
\item[$(i)$]  $A$ is an abelian subgroup of index $2$ in $G$;
\item[$(ii)$] $G$ is isomorphic to  the semidirect product $A \rtimes \Z / 2\Z$
where $\Z/ 2 \Z$ acts on $A$ by taking the inverse;
\item[$(iii)$] there is at most one element of order $2$ in $A$;
\item[$(iv)$] $A$ is a limit of cyclic groups.

\end{itemize}

Let us prove $(i)$. By $(P_1)$ of Lemma \ref{LemP}, $A$ is generated
by a set of pairwise commuting elements. Hence $A$ is abelian. As
$G$ is not abelian, the sentence $(P_1)$ implies that $G$ has at
least one non central element of order $2$. Let $s$ be such an
element. We show that $G=A \sqcup sA$. Let $x$ be in $G \setminus
A$. There are two cases:

\begin{itemize}
\item[(case 1)] $(sx)^2 \neq 1$. Then $sx$ belongs to $A$;
\item[(case 2)] $(sx)^2=1$. Since $s^2=x^2=1$, $s$ and $x$ commute. Hence $s$ and $x$ are non
central commuting elements of order $2$. We deduce from $(P_3)$ that
$sx$ belongs to $Z(G) \subset A$.
\end{itemize}
Thus $x$ belongs to $sA$ in both cases, which shows that $A$ has
index $2$ in $G$. As $G$ is finitely generated, so is $A$.

By $(P_2)$, central elements of $G$ have order at most $2$. The
sentence $(P_3)$ shows then that the conjugation by $s$ of an
element of $A$ consists in taking its inverse. Hence $(ii)$ is
proved.

Let us prove $(iii)$. Using $(P_4)$, we deduce that the center of
$G$ has at most two elements. We now show that elements of $A$ which
have order $2$ are central in $G$. Let $a$ be in $A$ and such that
$a^2=1$.  Write $a=yz$ with $y$ in the subgroup of $G$ generated by
$\{x \in G \, \vert \, x^2 \neq 1\}$ and $z$ in $Z(G)$. By $(P_2)$,
we have $s^{-1}as=y^{-1}z$. Since $a^2=y^2=1$, we deduce that
$s^{-1}as=a$. Thus $a$ is central, which completes the proof of
$(iii)$.

We now prove $(iv)$. We set $A^2:=\{a^2 \, \vert \, a \in A\}$ and
$\Dn^2:=\{g^2 \, \vert \, g \in \Dn\}$. We first show that $A^2$ is
a limit of cyclic groups. By proposition \ref{PropLimZn}, it
suffices to show that  $\Thu(A^2) \supset \bigcap_n \Thu(\Z / n\Z)$.
Let $\phi(x_1,\dots,x_k)$ be a quantifier free formula in variables
$x_1,\dots,x_k$ and consider the sentences $(P) \forall x_1\dots
\forall x_k \, \phi(x_1,\dots,x_k)$ and $ (P^2) \forall x_1\dots
\forall x_k \, \phi(x_1^2,\dots,x_k^2)$.
 We observe the following equivalences:
 $$\Dn \models P^2\Longleftrightarrow\Dn^2 \models P \mbox{ and }A \models P^2\Longleftrightarrow A^2
 \models P.$$Assume $\Z /n \Z \models P$ for all $n \ge 1$. Then $\Dn^2 \models
P$ for all $n \ge 3$ because $\Dn^2$ is a finite cyclic group. It
follows that $\Dn \models P^2$ for all $n \ge 3$. By assumption, $G
\models P^2$, hence $A \models P^2$. Consequently, $A^2 \models P$.
We deduce that $A^2$ is a limit of cyclic group, hence $A^2$ is
isomorphic to $\Z^n \times \Z / k \Z$ for some $n \ge 0, k \ge 1$ by
Proposition \ref{PropLimZn}. By $(P_4)$, there is at most one
element of order $2$ in $A$. We deduce that $A$ is isomorphic to
$\Z^n \times \Z / k'\Z$ with $k'$ in $\{k,2k\}$.
\end{sproof}

\section{The space of limits of dihedral groups on $m$ generators} \label{SecTop}

Let $G$ be a non abelian limit of dihedral groups that is generated
by two of its (necessarily non trivial) elements, say $x$ and $y$.
It follows from $(P_1)$ of Lemma \ref{LemP} that either $x$ or $y$
has order $2$. Assume then $x^2=1$. By $(P_3)$ of the same lemma, we
have either $y^2=1$ or simultaneously $y^2 \neq 1$ and
$x^{-1}yx=y^{-1}$. Thus $G$ is an homomorphic image of $\Di$ (both
$\Pres{x,y}{x^2=y^2=1}$ and $\Pres{x,y}{x^2=1,x^{-1}yx=y^{-1}}$ are
presentations of $\Di$). Hence $G$ is a dihedral group. This shows
that $2$-generated limits of dihedral groups are dihedral groups
(more generally, the rank of a limit $G=A \rtimes \Z/2\Z$ of
dihedral groups is $r(A)+1$). Consequently:

\begin{corollary}
The space of dihedral marked groups on $2$ generators is a closed
and open subspace of $\M$.
\end{corollary}

\begin{proof} It only remains to show that this subspace is open. As a quotient
of a dihedral group is a dihedral group, the result follows from
Lemma \ref{LemPresFinQuo}.
\end{proof}
 We have a complete topological description of
the space of dihedral groups on two generators. For each $n \neq 2$,
there are exactly three distinct marked groups on two generators
which are isomorphic to $\Dn$: $\begin{array}{lll}
A_{2n}&:=&\Pres{a,b}{ a^2=b^2=(ab)^n=1},\\ B_{2n}&:=&\Pres{a,b
}{a^2=b^n=1,a^{-1}ba=b^{-1}},\\
\ov{B}_{2n}&:=&\Pres{a,b}{b^2=a^n=1,b^{-1}ab=a^{-1}}.
\end{array}$\\ There is a unique marked group of $\M$ isomorphic to
$\mathbb{D}_4$: $A_4=B_4=\ov{B}_4$. The only accumulation points in
the space of dihedral groups on two generators are the three
distinct marked infinite dihedral groups $A_{\infty},B_{\infty}$ and
$\ov{B}_{\infty}$ (see Figure \ref{FigD2} below). This last fact is
proved in Proposition \ref{PropNbrCB}. The remaining statements can
be readily adapted from this proposition.

\begin{figure}[h]
\leavevmode \SetLabels
\L(.24*.5) $A_{\infty}$ \\
\L(.41*.53) $A_{2n}$ \\
\L(.55*.53) $A_4$ \\
\L(0.72*1) $B_{\infty}$ \\
\L(.64*.7) $B_{2n}$ \\
\L(0.64*0.28)  $\overline{B}_{2n}$\\
\L(0.72*0)  $\overline{B}_{\infty}$\\
\endSetLabels
%\ShowGrid
\begin{center}
\AffixLabels{\centerline{\epsfig{file =
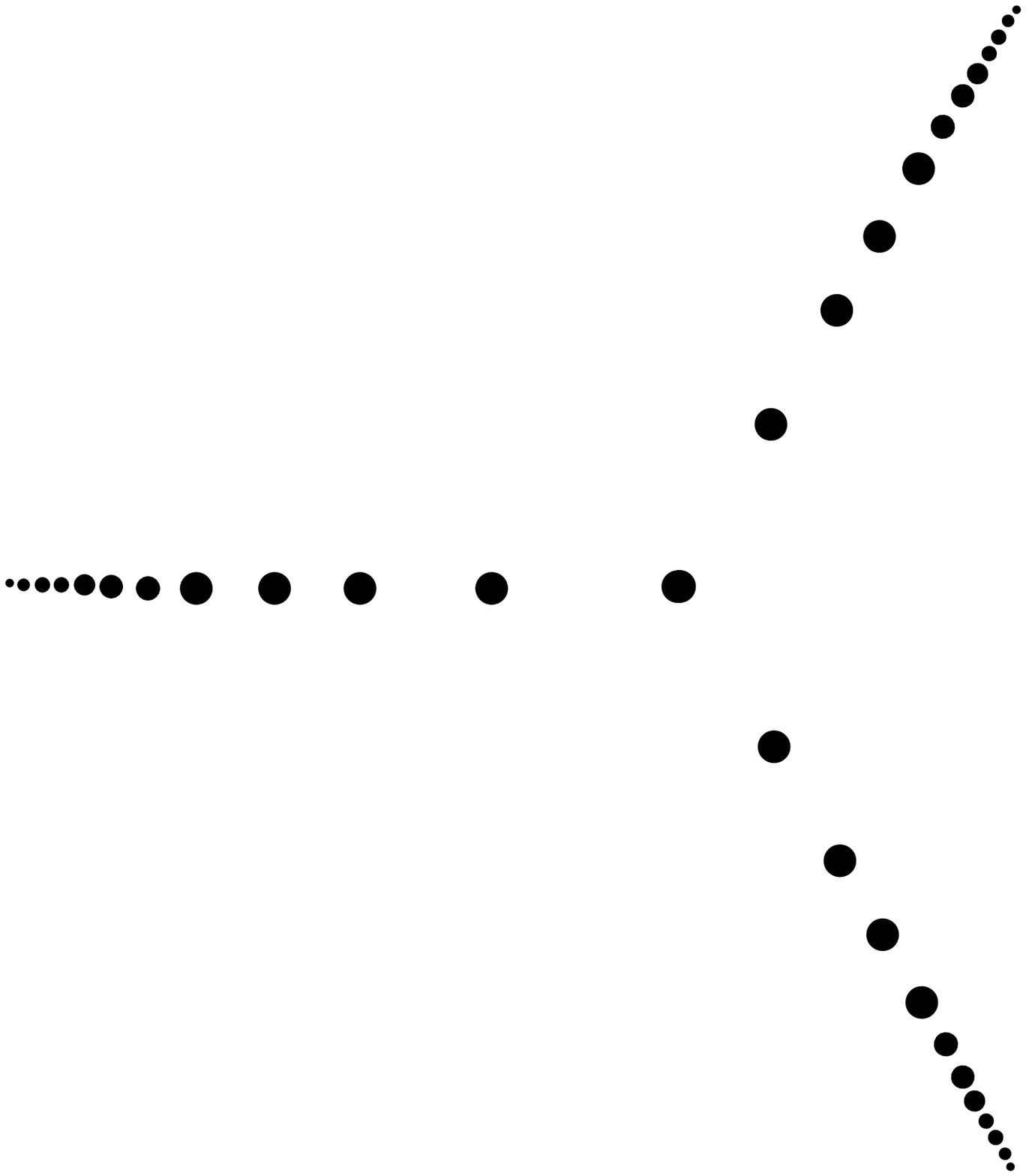,width=6cm,angle=0}}}
\end{center}
\caption{dihedral marked groups on $2$ generators} \label{FigD2}
\end{figure}
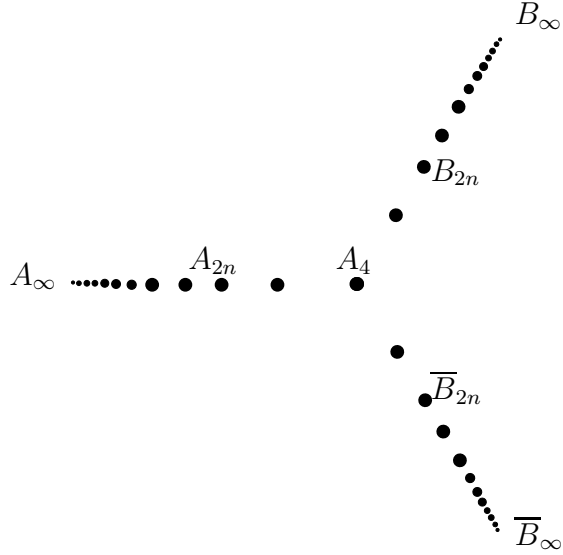

We can carry out such an analysis in the space of marked groups on
$m$ generators by using Cantor-Bendixson invariants defined in
Section \ref{SecCB}. We denote by $\omega$ the smallest infinite
ordinal, \it{i.e.} \rm the set $\N$ of positive integers endowed
with its natural order.

\begin{theorem}\label{ThmTopDn}
The topological closure $\Dihm$ of dihedral marked groups in $\Mm$
is homeomorphic to $\omega^{m-1}(2^m-1)+1$ endowed with the order
topology.
\end{theorem}
Because of the Mazurkiewicz-Sierpinski Theorem (Section
\ref{SecCB}), it suffices to show that the $(m-1)$-th derived set
$\mathcal{D}^{(m-1)}_m$ of $\Dihm$ contains $2^m-1$ points. This is
carried out with the two following propositions:

\begin{proposition} \label{PropMaxCB}
Let $G=A \rtimes \Z/2\Z$ be a limit of dihedral groups on $m$
generators. Then the Cantor-Bendixson rank of $G$ in $\Dihm$ is the
free rank of $A$.
\end{proposition}
Hence, the only remaining marked groups in $\mathcal{D}^{(m-1)}_m$
are marked groups abstractly isomorphic to $\Z^{m-1} \rtimes \Z /2
\Z$. We now count them:
\begin{proposition} \label{PropNbrCB}
Let $G=\Z^{m-1} \rtimes \Z/ 2\Z$ with $m \ge 2$.
 In $\Mm$, there are exactly $2^m-1$ marked groups which are abstractly
isomorphic to $G$.
\end{proposition}

We show Propositions \ref{PropMaxCB} and \ref{PropNbrCB}. Theorem
\ref{ThmTopDn} then directly follows.

Let $\mathcal{A} \subset \Mu$ be the set of abelian marked groups,
let $\mathcal{D} \subset \Mu$ be the set of dihedral marked groups
and let $\widetilde{\mathcal{D}} \subset \Mu$ be the set of
generalized dihedral marked groups. Let $\mathcal{C}$ be the
topological closure of cyclic marked groups in $\Mu$. We define
$Dih(A,S)$ with $S=(a_1,a_2,\dots)$ as the marked group
$(Dih(A),S')$ with $S'=(a,a_1,a_2,\dots)$ where $a$ denotes the
generator of $\Z/2\Z$.

\begin{lemma} \label{LemDih}
The map $Dih: \mathcal{A} \longrightarrow \tD$ is a continuous and
open embedding. Moreover, $Dih(\Cyc)=\D$.
\end{lemma}

\begin{sproof}{Proof of Lemma \ref{LemDih}}
We fix words $v_1,\dots,v_k,w_1,\dots,w_l$ in $\Free_{m+1}$ for
which the exponent sum of $e_1$ is zero. We then define the system
$$(\Sigma):\left\{
\begin{array}{c} v_1=1  ,\hdots ,v_k = 1 ,\\ w_1 \neq 1 ,\hdots
,w_l \neq 1 .\end{array} \right.$$ Let
$$D:=\Pres{e_1,\dots,e_{m+1}}{e_1^2=1,e_1e_ie_1^{-1}=e_i^{-1},
\br{e_i,e_j}=1, i=2,\dots,m+1}.$$ We reduce the words $v_i,w_j$ in
$D$ for $i=1,\dots,k$ and $j=1,\dots,l$ to get words without symbols
$e_1^{\pm 1}$. We then shift the indices on the left ($e_i$ becomes
$e_{i-1}$) to obtain words $v_i',w_j'$ in $\Fm$. We define
$(\Sigma')$ by replacing $v_i$ by $v_i'$ and $w_j$ by $w_j'$ in
$(\Sigma)$. We consider the elementary open sets $O_{\Sigma'}
\subset \mathcal{A}$ and $O_{\Sigma} \subset \mathcal{D}$. Let
$(A,S)$ be in $\mathcal{A}$. It is trivial to check that $Dih(A,S)
\in O_{\Sigma} \Longleftrightarrow (A,S) \in O_{\Sigma'}$. Hence
$Dih$ is a continuous and open map.

Assume $Dih(A,S)=Dih(B,T)$ with $(A,S),(B,T)$ in $\mathcal{A}$.
There is an isomorphism $\phi :Dih(A) \longrightarrow Dih(B)$ such
that $\phi \cdot S=T$. Observe that for any abelian group $C$, $C$
is the characteristic subgroup of $Dih(C)$ generated by the set $\{c
\in Dih(C) \, \vert \, c^2 \neq 1 \mbox{ or } c \mbox{ is central in
} Dih(C) \}$. This shows that $\phi$ induces an isomorphism from $A$
onto $B$. Hence $(A,S)=(B,T)$ which proves the injectivity of $Dih$.
\end{sproof}

\begin{lemma} \label{LemAbel}
Let $A$ be a finitely generated abelian group (respectively a limit
of cyclic groups). Then the Cantor-Bendixson rank of $A$ in $\Mu$
(respectively in $\mathcal{C}$) is the free rank of $A$.
\end{lemma}

\begin{sproof}{Proof of Lemma \ref{LemAbel}}
Consider a finitely generated abelian group $A$. Since $A$ is
finitely presented, the Cantor-Bendixson rank of $A$ in $\Mu$ is the
Cantor-Bendixson rank of $\{0\}$ in $\mathcal{N}(A)$ by Lemma
\ref{LemPresFinQuo}. The proof is an induction on the free rank $r$
of $A$. We first show that the Cantor-Bendixson rank of $\{0\}$ is
not less than $r$. If $r=0$, then $A$ is finite. It follows that
$\mathcal{N}(A)$ is a finite discrete space in which $\{0\}$ is
obviously isolated. Assume $r \ge 1$. Consider an infinite cyclic
subgroup $\langle z \rangle$ of $A$. For all $n$ in $\N$, the
Cantor-Bendixson rank of $\langle z^n \rangle $ in $\mathcal{N}(A)$
is the Cantor-Bendixson rank of $\{0\}$ in $\mathcal{N}(A/\langle
z^n \rangle)$ by Lemma \ref{LemPresFinQuo}. By the induction
hypothesis, this rank is at least $r-1$. Since $\langle z^n \rangle$
tends to $\{0\}$ in $\mathcal{N}(A)$ as $n$ tends to infinity, the
Cantor-Bendixson of $\{0\}$ is at least $r$. We now show that the
Cantor-Bendixson rank of $A$ is not greater than $r$. It is clear if
$r=0$. Assume $r \ge 1$. Consider the set $V$ of subgroups of $A$
whose intersection with $Tor(A)$ is trivial. Then $V$ is an open
neighborhood of $\{0\}$ in $\mathcal{N}(A)$. Let $B \neq \{0\}$ be
in $V$. By the induction hypothesis, the Cantor-Bendixson rank of
$\{0\}$ in $\mathcal{N}(A/B)$ is at most $r-1$. Since this is also
the Cantor-Bendixson rank of $B$ in $\mathcal{N}(A)$, the
Cantor-Bendixson rank of $\{0\}$ in $\mathcal{N}(A)$ is at most $r$.

If $A$ is in $\mathcal{C}$, we then consider the set
$\mathcal{N}_{\mathcal{C}}(A)$ of subgroups $B$ of $A$ such that
$A/B$ is a limit of cyclic groups. For any infinite cyclic factor
$\langle z \rangle$ of $A$, the subgroup $\langle z^n \rangle$
belongs to $\mathcal{N}_{\mathcal{C}}(A)$ for all $n \ge 1$ coprime
with $\vert Tor(A) \vert$. We can hence apply the reasoning above to
such an $A$.
\end{sproof}

\begin{sproof}{Proof of Proposition \ref{PropMaxCB}}
Let $G=A \rtimes \Z/2\Z$ be a limit of dihedral groups on $m$
generators. As $(G,S)$ is in the image of $Dih$ for a suitable
ordered generating set $S$, its Cantor-Bendixson rank in $\Dihm$ is
the Cantor-Bendixson rank of $A$ in $\mathcal{C}$ by Lemma
\ref{LemDih}. This is the free rank of $A$ by Lemma \ref{LemAbel}.
\end{sproof}

We fix $m \ge 2$ and $G=\Z^{m-1} \rtimes \Z /2 \Z$. We denote by $a$
the generating element of the subgroup $\Z/2\Z$ and we use the
additive notation in the normal subgroup $\Z^{m-1}$ of $G$. Let
$\phi$ be in $Aut(G)$. Since $\phi(a)$ has order $2$, we can write
$\phi(a)=v(\phi)a$ with $v(\phi)$ in $\Z^{m-1}$. We denote by
$\Z^{m-1} \rtimes GL_{n-1}(\Z)$ the semidirect product where
$GL_{m-1}(\Z)$ acts on $\Z^{m-1}$ in the standard way. We denote by
$(e_i)_{1 \le i \le n-1}$ the canonical basis of $\Z^{m-1}$.
\begin{lemma} \label{LemAutDih}
The map $$\begin{array}{ccc} \Phi: Aut(G) & \longrightarrow &
\Z^{m-1} \rtimes GL_{m-1}(\Z) \\
\phi&\longmapsto&(v(\phi),\phi_{\vert \Z^{m-1}})\\
\end{array}$$ is an isomorphism.
\end{lemma}

\begin{sproof}{Proof of Lemma \ref{LemAutDih}}
Since $\Z^{m-1}$ is the subgroup of $G$ generated by elements of
infinite order, $\Z^{m-1}$ is a characteristic subgroup. Hence
$\Phi$ is well defined. Checking that $\Phi$ is a homomorphism is
routine. Consider $v$ in $\Z^{m-1}$ and $f$ in $GL_{m-1}(\Z)$. The
group $G$ has the presentation
$\Pres{a,e_1,\dots,e_{m-1}}{a^2=1,ae_ia^{-1}=e_i^{-1},
\br{e_i,e_j}=1, i,j=1,\dots,m-1}$. We use then Von Dyck's Theorem to
show that there is unique automorphism $\phi$ of $G$ such that
$\phi(a)=va$ and $\phi_{\vert \Z^{m-1}}=f$.
\end{sproof}

\begin{sproof}{Proof of Proposition \ref{PropNbrCB}} The set of marked groups which
are isomorphic to $G$ in $\Mm$ corresponds bijectively to the the
set of $Aut(G)$-orbits of $V(G,m)$ of the diagonal action. For
$S=(g_1,\dots,g_m)$ in $V(G,m)$, we define $I(S)
\subset\{1,\dots,m\}$ as the set of indices $i$ satisfying
$g_i^2=1$. Let $P$ be a non empty subset of $\{1,\dots,m\}$. We
denote by $V(P)$, the set of all generating $G$-vectors $S$ of
length $m$ such that $I(S)=P$. Clearly, the sets $V(P)$ are pairwise
disjoint $Aut(G)$-invariant sets and there are $2^m-1$ such sets. We
now prove that the action of $Aut(G)$ is transitive on $V(P)$ for
all $P$. Starting with the generating vector
$S_0=(e_1,\dots,e_{m-1},a)$ and using elementary Nielsen
transformations, we can get a generating system $S'$ such that
$I(S')=P$ for any non empty $P$. Since the actions of the Nielsen
group $Aut(\Fm)$ on $V(G,m)$ commutes with the action of $Aut(G)$,
it suffices to prove that $Aut(G) \cdot S_0=V(\{m\})$. Let
$S=(f_1,\dots,f_{m-1},va)$ be in $V(\{m\})$. It is easy to check
that any word on the elements of $S$ that belongs to $\Z^{m-1}$ can
actually be generated by the elements $f_1,\dots,f_{m-1}$ only.
Consequently, $(f_1,\dots,f_{m-1})$ is a basis of $\Z^{m-1}$ and
there is $f$ in $GL_{m-1}(\Z)$ such that $f(e_i)=f_i$ for
$i=1,\dots,m-1$. By Lemma \ref{LemAutDih}, there is $\phi$ in
$Aut(G)$ such that $\phi(a)=va$ and $\phi_{\vert \Z^{m-1}}=f$. Thus
we have $\phi \cdot S_0=S$.
\end{sproof}

%%%%%%%%%%%%%%%%%%%%%%%%%%%%%%%%%%%%%%%%%%%%%%%%%%%%%%%%%%%%%%%%%%%%%%%%%%%%%%%%%%%%%%%%%%%%%%%%%%%%%%%%%%%%%%%%%%%%%%%%%%%%%%%%%
%                                              REFERENCES
%%%%%%%%%%%%%%%%%%%%%%%%%%%%%%%%%%%%%%%%%%%%%%%%%%%%%%%%%%%%%%%%%%%%%%%%%%%%%%%%%%%%%%%%%%%%%%%%%%%%%%%%%%%%%%%%%%%%%%%%%%%%%%%%%
\bibliographystyle{alpha}
\bibliography{BiblioLimDn}
%%%%%%%%%%%%%%%
\end{document}